\pdfoutput=1
\RequirePackage{ifpdf}
\ifpdf 
\documentclass[pdftex]{sigma}
\else
\documentclass{sigma}
\fi

\newtheorem{Theorem}{Theorem}
\newtheorem*{Theorem*}{Theorem}

\newtheorem{Proposition}[Theorem]{Proposition}
 { \theoremstyle{definition}

\newtheorem*{Remarks}{Remarks} }

\font\tensf=cmss10

\def\fl#1{\lfloor#1\rfloor}
\def\poq#1#2{(#1;q)_#2}
\def\al{\alpha}

\begin{document}


\renewcommand{\thefootnote}{}

\newcommand{\arXivNumber}{2102.02360}

\renewcommand{\PaperNumber}{089}

\FirstPageHeading

\ShortArticleName{Proof of Two Multivariate $q$-Binomial Sums Arising in Gromov--Witten Theory}

\ArticleName{Proof of Two Multivariate $\boldsymbol{q}$-Binomial Sums\\ Arising in Gromov--Witten Theory\footnote{This paper is a~contribution to the Special Issue on Basic Hypergeometric Series Associated with Root Systems and Applications in honor of Stephen C.~Milne's 75th birthday. The~full collection is available at \href{https://www.emis.de/journals/SIGMA/Milne.html}{https://www.emis.de/journals/SIGMA/Milne.html}}}

\Author{Christian KRATTENTHALER}

\AuthorNameForHeading{C.~Krattenthaler}

\Address{Fakult\"at f\"ur Mathematik, Universit\"at Wien,\\
Oskar-Morgenstern-Platz~1, A-1090 Vienna, Austria}
\Email{\href{mailto:Christian.Krattenthaler@univie.ac.at}{Christian.Krattenthaler@univie.ac.at}}
\URLaddress{\url{http://www.mat.univie.ac.at/~kratt/}}

\ArticleDates{Received February 03, 2024, in final form October 07, 2024; Published online October 10, 2024}

\Abstract{We prove two multivariate $q$-binomial identities conjectured by Bousseau, Brini and van Garrel [{\it Geom.\ Topol.} {\bf 28} (2024), 393--496, arXiv:2011.08830] which give generating series for Gromov--Witten invariants of two specific log Calabi--Yau surfaces. The key identity in all the proofs is Jackson's $q$-analogue of the Pfaff--Saalsch\"utz summation formula from the theory of basic hypergeometric series.}


\Classification{33D15; 05A30; 14J32; 14N35; 53D45; 57M27}

\begin{flushright}
\begin{minipage}{53mm}
\it Dedicated to Stephen Milne,\\
 the pioneer of multidimensional\\
  basic hypergeometric series
\end{minipage}
\end{flushright}

\medskip

\renewcommand{\thefootnote}{\arabic{footnote}}
\setcounter{footnote}{0}

The purpose of this note is to prove two conjectured multivariate
$q$-binomial summation identities from \cite{BoBGAA}. There,
Bousseau, Brini and van Garrel deal with the computation of
Gromov--Witten invariants of log Calabi--Yau surfaces (Looijenga pairs).
For the two non-tame (but quasi-tame) surfaces $\text{dP}_1(0,4)$
and $\mathbb F_0(0,4)$, conjectured closed-form expressions are given
in \cite{BoBGAA} for the
corresponding generating series
(I refer to \cite{BoBGAA} for background and notation), namely
(cf.~\mbox{\cite[Conjecture~B.3, equation~(B-2)]{BoBGAA}})
\begin{equation} \label{AA}
\text{\tensf N}_{(d_0,d_1)}^{\log}(\text{dP}_1(0,4))(\hbar)
=
\frac {[2d_0]_q} {[d_0]_q}
\begin{bmatrix} d_0\\d_1\end{bmatrix}_q
\begin{bmatrix} d_0+d_1-1\\d_1-1\end{bmatrix}_q
\end{equation}
and (cf.~\cite[Conjecture~B.3, equation~(B-3)]{BoBGAA})
\begin{equation} \label{AB}
\text{\tensf N}_{(d_1,d_2)}^{\log}(\mathbb F_0(0,4))(\hbar)
=
\frac {[2d_1+d_2]_q} {[d_2]_q}\begin{bmatrix} d_1+d_2-1\\d_2-1\end{bmatrix}_q^2,
\end{equation}
where $q={\rm e}^{i\hbar}$.
Here, the $q$-integers $[\al]_q$ are defined symmetrically according
to physics convention, $[\al]_q:=q^{\al/2}-q^{-\al/2}$, and
the (corresponding) $q$-binomial coefficients
$\bigl[\smallmatrix n\\k\endsmallmatrix\bigr]_q$ are defined~by\footnote{Since this definition is based on the
$q$-integers according to physics convention,
these ``physics $q$-binomial coefficients'' differ
from the ``combinatorics $q$-binomial coefficients'' (cf.\
\cite[Exercise~1.2]{GaRaAF})
by a multiplicative factor of~${q^{-k(n-k)/2}}$.}
\begin{equation*}
\begin{bmatrix} n\\k\end{bmatrix}_q:=\begin{cases} \displaystyle
\frac {[n]_q[n-1]_q\cdots[n-k+1]_q} {[k]_q[k-1]_q\cdots[1]_q}
&\text{if }0\le k\le n,\\
0&\text{otherwise.}\end{cases}
\end{equation*}

In relation to the above two conjectures, we prove the following
two multivariate summation identities.

\begin{Theorem} \label{TA}
For integers $d_0$ and $d_1$ with $d_0>d_1\ge1$, we have
\begin{gather}
\sum_{m\ge1}
\underset{n_1>\dots>n_m>0}
{
\underset{k_1,\dots,k_m>0,\,k_0\ge0,}
{
\underset{n_1k_1+\dots+n_mk_m=d_0-d_1,}
{\sum_{k_1+\dots+k_m=d_1-k_0,}}}}
\begin{bmatrix} 2d_0\\k_1\end{bmatrix}_q
\begin{bmatrix} 2d_0-2(n_1-n_2)k_1\\k_2\end{bmatrix}_q\cdots\nonumber\\
\qquad{}\times
\begin{bmatrix} 2d_0-2\sum_{j=1}^{m-1}(n_j-n_m)k_j\\k_m\end{bmatrix}_q
\begin{bmatrix} 2d_1\\k_0\end{bmatrix}_q
=\frac {[2d_0]_q} {[d_0]_q}
\begin{bmatrix} d_0\\d_1\end{bmatrix}_q
\begin{bmatrix} d_0+d_1-1\\d_0\end{bmatrix}_q.
\label{AF}
\end{gather}
\end{Theorem}

\begin{Theorem} \label{TB}
For positive integers $d_1$ and arbitrary $d_2$, we have
\begin{gather}
\sum_{m\ge1}
\underset{n_1>\dots>n_m>0}
{
\underset{k_1,\dots,k_m>0,}
{\sum_{n_1k_1+\dots+n_mk_m=d_1,}}}
\begin{bmatrix} d_2+2d_1\\k_1\end{bmatrix}_q
\cdots
\begin{bmatrix} d_2+2d_1+2\sum_{j=1}^{i}(n_{i}-n_{j})k_{j}\\k_{i}\end{bmatrix}_q \cdots\nonumber\\
\qquad{}\times
\begin{bmatrix} d_2+2d_1+2\sum_{j=1}^{m}(n_{m}-n_{j})k_{j}\\k_{m}\end{bmatrix}_q
\begin{bmatrix} d_2\\\sum_{j=1}^mk_{j}\end{bmatrix}_q
\!=
\frac {[2d_1+d_2]_q} {[d_2]_q}\begin{bmatrix} d_1+d_2-1\\d_1\end{bmatrix}_q^2.
\label{AG}
\end{gather}
\end{Theorem}

\begin{Remarks}\quad
\begin{enumerate}\itemsep=0pt
\item[(1)] Both identities actually hold when $q$ is considered as a formal variable.
Furthermore, in the statement of Theorem~\ref{TB}, the phrase ``arbitrary~$d_2$"
means that the identity holds when~$d_2$~is considered as a formal variable.

\item[(2)] For any fixed $d_0$ and $d_1$,
the sums over $m$ in \eqref{AF} and \eqref{AG} are {\it finite} sums since
all of the $k_i$'s and the $n_i$'s are at least~1, which implies
an obvious bound on~$m$.
\end{enumerate}
\end{Remarks}

By \cite[Theorem~B.1]{BoBGAA}, Theorem~\ref{TA} implies~\eqref{AA}.
Similarly, by \cite[Theorem~B.2]{BoBGAA}, Theorem~\ref{TB}
implies~\eqref{AB}.\footnote{For the sake of consistency, in
comparison to Theorem~B.2 in \cite{BoBGAA},
here we have reversed the indexing of the~$k_i$'s and the~$n_i$'s,
that is, we have replaced $k_i$ by $k_{m-i+1}$ and $n_i$ by
$n_{m-i+1}$, $i=1,2,\dots,m$.}

The identities~\eqref{AA} and~\eqref{AB} are actually special cases of a more
general conjecture, namely \mbox{\cite[Conjecture~4.7]{BoBGAA}}, which
predicts a closed-form formula for
$\text{\tensf
N}_{(d_0,d_1,d_2,d_3)}^{\log}(\text{dP}_3(0,2))(\hbar)$.
It~is~con\-ceiv\-able that the ideas of this note, or similar ones,
may lead to a proof of this more general conjecture. However, as is
explained in \cite[Section~4.2]{BoBGAA}, in order to obtain an
expression for~$\text{\tensf
N}_{(d_0,d_1,d_2,d_3)}^{\log}(\text{dP}_3(0,2))(\hbar)$
to start with,
one would have to perform certain scattering diagram calculations.
This is deemed ``daunting" by the authors of~\cite{BoBGAA} (see the
paragraph below Conjecture~4.7), and they do not carry out these
calculations.

As is the case frequently, the identities in Theorems~\ref{TA} and~\ref{TB}
are difficult (impossible?) to prove directly since the parameters
in these identities do not allow for enough flexibility, in particular
if one has an inductive approach in mind (which we do).
The key in proving~\eqref{AF}
and~\eqref{AG} is to {\it generalise}, or, in this case, to {\it refine}.
By experimenting with the sums in~\eqref{AF} and~\eqref{AG}, I noticed that
one can still get closed forms if we fix the sum of the~$k_i$'s,
$i=1,2,\dots,m$. This leads us to the following key result.

\begin{Proposition} \label{TC}
Let $k_0$ and $d_1$ be integers with $1\le k_0\le d_1$.
Furthermore, for arbitrary~$d_0$ set
\begin{equation} \label{AC}
f(d_0,d_1,k_0)=\sum_{m\ge1}
\underset{n_1>\dots>n_m>0}
{
\underset{k_1,\dots,k_m>0,}
{
\underset{n_1k_1+\dots+n_mk_m=d_1,}
{\sum_{k_1+\dots+k_m=k_0,}}}}
\prod_{i=1}^m \begin{bmatrix}
  2d_0-2\sum_{j=1}^{i-1}(n_j-n_i)k_j\\k_i\end{bmatrix}_q .
\end{equation}
Then
\begin{equation*}
f(d_0,d_1,k_0)=\frac {[2d_0]_q} {[k_0]_q}
\begin{bmatrix} 2d_0-d_1+k_0-1\\k_0-1\end{bmatrix}_q
\begin{bmatrix} d_1-1\\
k_0-1\end{bmatrix}_q.
\end{equation*}
\end{Proposition}

\begin{Remarks}\quad
\begin{enumerate}\itemsep=0pt
\item[(1)] Again, the identity actually holds when $q$ is considered as a
formal variable.

\item[(2)]
For the meaning of ``arbitrary $d_0$", see Remarks~(1) below
Theorem~\ref{TB}.
\end{enumerate}
\end{Remarks}

Before we can embark on the proof of the proposition, we need to
introduce the standard notation for basic hypergeometric series,
\begin{equation} \label{AH}
{}_{r+1}\phi_r\!\left[\begin{matrix} a_1,\dots,a_{r+1}\\
b_1,\dots,b_r\end{matrix}; q,z\right]
=\sum _{\ell=0} ^{\infty}\frac {\poq{a_1}{\ell}\cdots\poq{a_{r+1}}{\ell}}
{\poq{q}{\ell}\poq{b_1}{\ell}\cdots\poq{b_r}{\ell}}
z^\ell,
\end{equation}
where $(a;q)_0=1$ and $(a;q)_m=\prod_{k=0}^{m-1}\bigl(1-aq^k\bigr)$.
The ``bible" \cite{GaRaAF} of the theory of basic hypergeometric
series contains many summation and transformation formulae for such
series. The formula that we need here is Jackson's $q$-analogue
of the Pfaff--Saalsch\"utz summation
(see~\mbox{\cite[equation~(1.7.2); Appendix II.12]{GaRaAF}})
\begin{equation} \label{AE}
{}_3\phi _2\!\left [ \begin{matrix} \let\over/ a,b,{q^{-N}}\\ \let\over/  c,{{a b {q^{1 - N}}}/
   c}\end{matrix} ;q,q\right ] ={\frac {{(\let\over/ {c/ a};q)}_{N}\, {(\let\over/ {c/ b};q)}_{N}}
    {{(\let\over/ c;q)}_{N}\, {(\let\over/ {c\over {a b}};q)}_{N}}},
\end{equation}
where $N$ is a nonnegative integer.

\begin{proof}[Proof of Proposition \ref{TC}]
We prove the claim by induction on $k_0$.

First we consider the start of the induction, $k_0=1$.
In this case, the summation on  the right-hand side of \eqref{AC} reduces to
$m=1$,  $k_1=1$, $n_1=d_1$, and hence
\begin{equation*}
f(d_0,d_1,1)=\begin{bmatrix} 2d_0\\1\end{bmatrix}_q,
\end{equation*}
in agreement with our assertion.

For the induction step, we rewrite the definition of $f(d_0,d_1,k_0)$
in~\eqref{AC} in the form
\begin{gather}
f(d_0,d_1,k_0)=
\chi(k_0\mid d_1)\begin{bmatrix} 2d_0\\k_0\end{bmatrix}_q
\nonumber\\ \qquad{}
+\sum_{m\ge2}\sum_{k_m=1}^{k_0-1}\sum_{n_m=1}^{\fl{(d_1-k_0+k_m)/k_0}}
\kern-3pt
\underset{\bar n_1>\dots>\bar n_{m-1}>0}
{
\underset{k_1,\dots,k_{m-1}>0,}
{
\underset{\bar n_1k_1+\dots+\bar n_{m-1}k_{m-1}=d_1-n_mk_0,}
{\sum_{k_1+\dots+k_{m-1}=k_0-k_m,}}}}
\kern-3pt
\begin{bmatrix} 2d_0\\k_1\end{bmatrix}_q
\begin{bmatrix} 2d_0-2(\bar n_1-\bar n_2)k_1\\k_2\end{bmatrix}_q
\cdots\nonumber\\ \qquad\quad{}
\times
\begin{bmatrix} 2d_0-2\sum_{j=1}^{m-2}(\bar n_j-\bar n_{m-1})k_j\\k_{m-1}\end{bmatrix}_q
\begin{bmatrix} 2d_0-2d_1+2n_mk_0\\k_m\end{bmatrix}_q,
\label{AEa}
\end{gather}
where $\bar n_i=n_i-n_m$, $i=1,2,\dots,m-1$, and where $\chi(\cdot)$
denotes the usual truth function, that is, $\chi(\mathcal A)=1$ if
$\mathcal A$ is true and $\chi(\mathcal A)=0$ otherwise.
There are two details in this expression which require further
explanation. First of all, the first term on the right-hand side
of~\eqref{AEa} gives the contribution of the sum in~\eqref{AC} for $m=1$.
Indeed, for $m=1$ the summation conditions in~\eqref{AC} require
$k_1=k_0$ and $n_1k_1=d_1$. Consequently, a corresponding summand
occurs only if $k_0=k_1\mid d_1$; if this divisibility condition
is satisfied, the summand equals $\bigl[\smallmatrix
2d_0\\k_0\endsmallmatrix\bigr]_q$. Second, by~the~conditions imposed
on $k_1,\dots,k_{m-1}$ and $\bar n_1,\dots,\bar n_{m-1}$, the sum
$\bar n_1k_1+\dots+\bar n_{m-1}k_{m-1}$ must be at least as large
as the sum $k_1+\dots+k_{m-1}$. This explains the upper bound on~$n_m$
in the inner sum.

Now, again with \eqref{AC} in mind, the equation \eqref{AEa} may be written as
\begin{gather*}
f(d_0,d_1,k_0)=
\chi(k_0\mid d_1)\begin{bmatrix} 2d_0\\k_0\end{bmatrix}_q\\ \hphantom{f(d_0,d_1,k_0)=}{}
+\sum_{k=1}^{k_0-1}
\sum_{n=1}^{\fl{(d_1-k_0+k)/k_0}}
f(d_0,d_1-nk_0,k_0-k)
\begin{bmatrix} 2d_0-2d_1+2nk_0\\k\end{bmatrix}_q.
\end{gather*}
We may now use the induction hypothesis, and obtain
\begin{gather}
\notag
f(d_0,d_1,k_0)=\chi(k_0\mid d_1)\begin{bmatrix} 2d_0\\k_0\end{bmatrix}_q
+\sum_{k=1}^{k_0-1}
\sum_{n=1}^{\fl{(d_1-k_0+k)/k_0}}
\frac {[2d_0]_q} {[k_0-k]_q} \\ \notag
\hphantom{f(d_0,d_1,k_0)=}{}
\times\begin{bmatrix} 2d_0-d_1+nk_0+k_0-k-1\\k_0-k-1\end{bmatrix}_q
\begin{bmatrix} d_1-nk_0-1\\
k_0-k-1\end{bmatrix}_q
\begin{bmatrix} 2d_0-2d_1+2nk_0\\k\end{bmatrix}_q\\ \notag \hphantom{f(d_0,d_1,k_0)}{}
=\chi(k_0\mid d_1)\begin{bmatrix} 2d_0\\k_0\end{bmatrix}_q
+\sum_{n=1}^{\fl{(d_1-1)/k_0}}
\sum_{k=0}^{k_0-1}
\frac {[2d_0]_q} {[k_0-k]_q} \\ \notag \hphantom{f(d_0,d_1,k_0)=}{}
\times\begin{bmatrix} 2d_0-d_1+nk_0+k_0-k-1\\k_0-k-1\end{bmatrix}_q
\begin{bmatrix} d_1-nk_0-1\\
k_0-k-1\end{bmatrix}_q
\begin{bmatrix} 2d_0-2d_1+2nk_0\\k\end{bmatrix}_q\\ \hphantom{f(d_0,d_1,k_0)=}{}
-\sum_{n=1}^{\fl{(d_1-1)/k_0}}
\frac {[2d_0]_q} {[k_0]_q}
\begin{bmatrix} 2d_0-d_1+nk_0+k_0-1\\k_0-1\end{bmatrix}_q
\begin{bmatrix} d_1-nk_0-1\\
k_0-1\end{bmatrix}_q.
\label{AD}
\end{gather}
Now we write the sum over~$k$
in terms of the standard basic hypergeometric notation~\eqref{AH}.
Thus, this sum over $k$ becomes
\begin{gather*}
\frac {[2d_0]_q} {[k_0]_q}
\begin{bmatrix} 2d_0-d_1+nk_0+k_0-1\\k_0-1\end{bmatrix}_q
\begin{bmatrix} d_1-nk_0-1\\
k_0-1\end{bmatrix}_q
\\ \qquad{}\times
   {}_3\phi _2 \!\left[\begin{matrix} q^{-2 d_0 + 2 d_1 - 2 k_0 n},
q^{-k_0},q^{1-k_0}\\
q^{1 + d_1 - k_0 - k_0 n},q^{1 - 2 d_0 + d_1 - k_0 - k_0 n}
\end{matrix} ;q,q\right].
\end{gather*}
The $_3\phi_2$-series can be evaluated by means of the
$q$-Pfaff--Saalsch\"utz summation~\eqref{AE}.
After simpli\-fi\-ca\-tion, we arrive at the expression
\begin{equation*}
\frac {[2d_0]_q} {[k_0]_q}
\begin{bmatrix} 2d_0-d_1+k_0n-1\\k_0-1\end{bmatrix}_q
\begin{bmatrix} d_1-k_0n+k_0-1\\k_0-1\end{bmatrix}_q.
\end{equation*}
If we substitute this in \eqref{AD}, then we get
\begin{gather*}
f(d_0,d_1,k_0)=
\chi(k_0\mid d_1)\begin{bmatrix} 2d_0\\k_0\end{bmatrix}_q\\ \hphantom{f(d_0,d_1,k_0)=}{}
+\sum_{n=1}^{\fl{(d_1-1)/k_0}}
\frac {[2d_0]_q} {[k_0]_q}
\begin{bmatrix} 2d_0-d_1+k_0n-1\\k_0-1\end{bmatrix}_q
\begin{bmatrix} d_1-k_0n+k_0-1\\k_0-1\end{bmatrix}_q\\ \hphantom{f(d_0,d_1,k_0)=}{}
-\sum_{n=1}^{\fl{(d_1-1)/k_0}}
\frac {[2d_0]_q} {[k_0]_q}
\begin{bmatrix} 2d_0-d_1+nk_0+k_0-1\\k_0-1\end{bmatrix}_q
\begin{bmatrix} d_1-nk_0-1\\
k_0-1\end{bmatrix}_q.
\end{gather*}
As is straightforward to verify, the first term in this expression
can be integrated into the first sum, so that we obtain
\begin{gather}
f(d_0,d_1,k_0)=
\sum_{n=1}^{\fl{d_1/k_0}}
\frac {[2d_0]_q} {[k_0]_q}
\begin{bmatrix} 2d_0-d_1+k_0n-1\\k_0-1\end{bmatrix}_q
\begin{bmatrix} d_1-k_0n+k_0-1\\k_0-1\end{bmatrix}_q \nonumber \\ \hphantom{f(d_0,d_1,k_0)=}{}
-\sum_{n=1}^{\fl{(d_1-1)/k_0}}
\frac {[2d_0]_q} {[k_0]_q}
\begin{bmatrix} 2d_0-d_1+nk_0+k_0-1\\k_0-1\end{bmatrix}_q
\begin{bmatrix} d_1-nk_0-1\\
k_0-1\end{bmatrix}_q.
\label{ADa}
\end{gather}
On the right-hand side,
these are, up to a shift of the index~$n$ and the slightly deviating
lower and upper bounds on the summation index, essentially the same
sums. After cancellation, only the term corresponding to $n=1$ of
the first sum and the term corresponding to $n=\fl{(d_1-1)/k_0}$
of the second sum survives, however the latter only if $k_0\nmid d_1$.
(If $k_0\mid d_1$, then no term of the second sum survives the
cancellation.)
Since, for $n=\fl{(d_1-1)/k_0}$, we have
\begin{equation*}0\le d_1-nk_0-1=d_1-\biggl\lfloor{\frac {d_1-1} {k_0}}\biggr\rfloor k_0-1
<d_1-\biggl(\frac {d_1} {k_0}-1\biggr)k_0-1
= k_0-1,\end{equation*}
the $q$-binomial coefficient $\bigl[\smallmatrix d_1-nk_0-1\\
k_0-1\endsmallmatrix\bigr]_q$ vanishes for this choice of parameters.
In other words, after cancellation of terms in~\eqref{ADa}, only the
term corresponding to $n=1$ of the first sum survives and gives a non-zero
contribution. These arguments yield
\begin{equation*}
f(d_0,d_1,k_0)=
\frac {[2d_0]_q} {[k_0]_q}
\begin{bmatrix} 2d_0-d_1+k_0-1\\k_0-1\end{bmatrix}_q
\begin{bmatrix} d_1-1\\k_0-1\end{bmatrix}_q.
\end{equation*}
This completes the induction step and the proof of the
proposition.
\end{proof}

Now we are in the position to prove Theorems~\ref{TA} and~\ref{TB}.

\begin{proof}[Proof of Theorem \ref{TA}]
The left-hand side of \eqref{AF} can be
rewritten as
\begin{equation*}
\sum_{k_0=0}^{d_1}f(d_0,d_0-d_1,d_1-k_0)\begin{bmatrix} 2d_1\\k_0\end{bmatrix}_q.
\end{equation*}
If we now use Proposition~\ref{TC} for the evaluation of
$f(d_0,d_0-d_1,d_1-k_0)$ and write the result
in basic hypergeometric notation \eqref{AH}, the above sum becomes
\begin{equation*}
\frac {[2d_0]_q} {[d_1]_q}
\begin{bmatrix} d_0+2d_1-1\\d_1-1\end{bmatrix}_q
\begin{bmatrix} d_0-d_1-1\\
d_1-1\end{bmatrix}_q
   {}_3\phi _2\!\left [ \begin{matrix} q^{-2d_1},q^{-d_1},q^{1-d_1}\\
q^{1+d_0-2d_1},q^{1-d_0-2d_1}
\end{matrix} ;q,q\right ].
\end{equation*}
Also this $_3\phi_2$-series can be evaluated by means of the
$q$-Pfaff--Saalsch\"utz summation~\eqref{AE}.
After simplification, one obtains
\begin{equation*}
\frac {[2d_0]_q} {[d_0]_q}
\begin{bmatrix} d_0\\d_1\end{bmatrix}_q
\begin{bmatrix} d_0+d_1-1\\d_0\end{bmatrix}_q,
\end{equation*}
as desired.
\end{proof}

\begin{proof}[Proof of Theorem \ref{TB}]
With $k_0=\sum_{j=1}^mk_j$,
the sum on the left-hand side of~\eqref{AG} can be rewritten as
\begin{equation*}
\sum_{k_0\ge1}f\biggl(d_1+\frac {d_2} {2},d_1,k_0\biggr)
\begin{bmatrix} d_2\\k_0\end{bmatrix}_q.
\end{equation*}
If we now use Proposition~\ref{TC} for the evaluation of
$f(d_1+\tfrac {d_2} {2},d_1,k_0)$ and write the result
in basic hypergeometric notation \eqref{AH}, the above sum becomes
\begin{equation*}
\frac {[2d_1+d_2]_q\,[d_2]_q} {[1]_q^2}
   {}_3\phi _2\!\left [ \begin{matrix} q^{1+d_1+d_2},q^{1-d_2},q^{1-d_1}\\
q^{2},q^{2}
\end{matrix} ;q,q\right ].
\end{equation*}
Again, the $_3\phi_2$-series can be evaluated by means of the
$q$-Pfaff--Saalsch\"utz summation~\eqref{AE}. As~a~result, we obtain
\begin{equation*}
\frac {[2d_1+d_2]_q} {[d_2]_q}\begin{bmatrix} d_1+d_2-1\\d_1\end{bmatrix}_q^2,
\end{equation*}
as desired.
\end{proof}

We close with the remark that it is highly surprising that in all
three proofs the identity from the theory of basic hypergeometric
series that is required is the $q$-Pfaff--Saalsch\"utz
summation~\eqref{AE}. Usually, one needs the $q$-Chu--Vandermonde summation
here, a transformation formula there, and maybe the $q$-Pfaff--Saalsch\"utz
summation somewhere. However, remarkably, here it is {\it exclusively}
the $q$-Pfaff--Saalsch\"utz summation.

\subsection*{Acknowledgements}

I am indebted to the anonymous referees for their many extremely
helpful comments, which not only led to an improved presentation
but also let me discover a subtle gap in the proof of Proposition~\ref{TC}
in the original version of the manuscript.
This research was partially supported by the Austrian
Science Foundation FWF (grant S50-N15)
in the framework of the Special Research Program
``Algorithmic and Enumerative Combinatorics''.

\pdfbookmark[1]{References}{ref}
\LastPageEnding

\end{document}